# Magic Knight's Tours in Higher Dimensions

Awani Kumar. 2/11-C, Vijayant Khand, Gomti Nagar, Lucknow 226 010. INDIA
E-mail: awanieva@gmail.com

**Abstract**

A knight's tour on a board is a sequence of knight moves that visits each square exactly once. A knight's tour on a square board is called magic knight's tour if the sum of the numbers in each row and column is the same (magic constant). Knight's tour in higher dimensions ($n > 3$) is a new topic in the age-old world of knight's tours. In this paper, it has been proved that there can't be magic knight's tour or closed knight's tour in an odd order $n$-dimensional hypercube. 3 x 4 x $2^{n-2}$ is the smallest cuboid ($n \geq 2$) and 4 x 4 x $4^{n-2}$ is the smallest cube in which knight's tour is possible in $n$-dimensions ($n \geq 3$). Magic knight's tours are possible in 4 x 4 x 4 x 4 and 4 x 4 x 4 x 4 x 4 hypercube.

## 1 Introduction

"The oldest of knight puzzles is the knight's tour," asserts Martin Gardner [1]. The problem is more than 1000 years old. The chess historian H. J. R. Murray [2] describes closed tours of the 8×8 board by the Shatranj players al-Adli and as-Suli who lived in Baghdad around 840 and 900 CE respectively. Knight's tour questions have continued to fascinate both amateur and professional mathematicians ever since. The mathematician A.T. Vandermonde [3] was the first to construct a three-dimensional knight's tour, in a 4 × 4 × 4 cube, published in 1771. Other 3D examples have been provided by Schubert [4], Gibbins [5], Stewart [6], Jelliss [7], Petkovic [8] and DeMaio [9]. More recently the present author, Awani Kumar [10] [11], looked into the possibilities of knight's tours in cubes and cuboids, having "magic" properties (i.e. adding to a fixed sum along most lines); some new results are reported below. The natural next step is to extend the study to knight's tours in higher dimensions. Is it possible to find magic knight tours in a 4D hypercube? Can such magic tours be extended further into five (or still higher) dimensions? Is it possible to have magic knight's tour or closed knight's tour in an odd order hypercube? What is the size of the smallest cuboid and smallest cube in which knight's tour is possible in $n$-dimensions? The author plans to answer these original questions. We describe such a higher dimensional construction, divided into lines of $m$ cells, as an $n$-dimensional "lattice" of "order" $m$.

## 2 Knight's tours in three dimensions

Knight moves two squares horizontally and one square vertically, or one square horizontally and two squares vertically on a plane board (that is, in two dimensions). It is shown in Figure 1a. A *knight's tour* is a sequence of knight moves that visits each square exactly once. Figure 1b shows knight's tour on 6 x 6 board. The discerning readers must have observed that, here, all the rows add up to 111. A tour of a knight on a square board is called *magic knight's tour* if the sum of the numbers in each row and column is the same (the *magic constant*). In three dimensions the knight is assumed to move in its usual fashion in each of the three mutually perpendicular planes through the cell it initially occupies. If the three coordinate directions are *x*, *y*, *z* then the three planes can be represented by the pairs of coordinates *xy*, *xz*, *yz*. Thus the mobility of the knight is multiplied three-fold; on a two-dimensional board the knight has a choice of up to 8 cells to which it can move, but in three dimensions it can have as many as 24 cells to move. However, on small boards, or near the edges of larger boards, the number of moves will of course be less than this maximum, since blocked by the board edges or faces. The knight cannot move at all in a 2 × 2 × 2 cube, or from the central cell of a 3 × 3 × 3 cube. So the smallest cubical board on which the knight is mobile on every cell is the 4 × 4 × 4 cube.

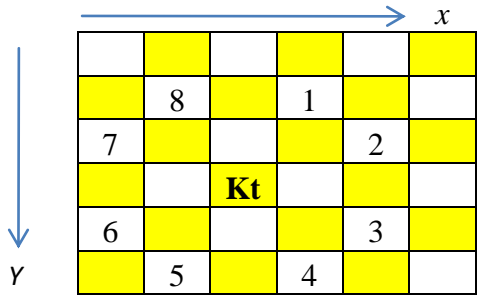

| 1  | 10 | 35 | 28 | 25 | 12 |
|----|----|----|----|----|----|
| 34 | 29 | 2  | 11 | 8  | 27 |
| 3  | 36 | 9  | 26 | 13 | 24 |
| 30 | 33 | 20 | 5  | 16 | 7  |
| 21 | 4  | 31 | 18 | 23 | 14 |
| 32 | 19 | 22 | 15 | 6  | 17 |

Figure 1a. Possible knight's move on a plane board    Figure 1b. Knight's tour on 6 x 6 board

Figure 1c shows the 12 possible moves of a knight from a central cell of a 4 × 4 × 4 cube. Readers can visualize this 3D board as a stack of four 2D layers, one above the other, lettered in alphabetical order A to D. As in the 2D case the cubical cells (or their floors) can be coloured alternately light and dark so that a knight at a light cell can only jump to a dark cell and vice versa. This remains true in higher dimensions. Figure 2 is an example of knight's tour in a 4 × 4 × 4 cube. Here, the sum of all the rows, columns and pillars are the same (equal to the magic constant 130). This magic knight's tour is closed, since its initial cell (1) and last cell (64) are connected by a knight's move; on the other hand Figure 3 is an open magic tour. It will be found that the sums of the space diagonals of these cubes are also equal to the magic constant 130, so they are "magic cubes" in the traditional sense. Stertenbrink [12] was the first to discover magic tour of knight in a 4 × 4 × 4 cube in the year 2003. Subsequently Aale de Winkel [13] has compiled 142 such magic tours in an attempt to enumerate all the distinct solutions. The author has found larger magic knight tours by extending the pattern to 8 × 8 × 8 and 12 × 12 × 12 cubes.

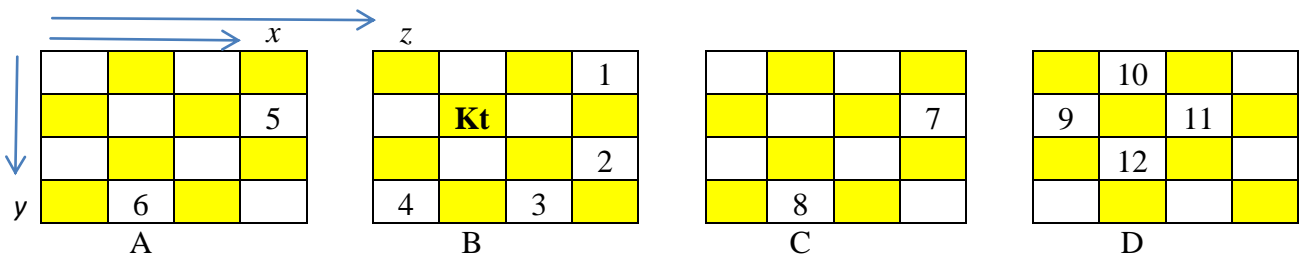

Figure 1c. Possible knight's move in 4 x 4 x 4 cube

| 19 | 46 | 63 | 2 |
|----|----|----|----|
| 48 | 1  | 20 | 61 |
| 29 | 52 | 33 | 16 |
| 34 | 31 | 14 | 51 |

A

| 58 | 7  | 22 | 43 |
|----|----|----|----|
| 21 | 60 | 41 | 8  |
| 40 | 9  | 28 | 53 |
| 11 | 54 | 39 | 26 |

B

| 47 | 18 | 3  | 62 |
|----|----|----|----|
| 4  | 45 | 64 | 17 |
| 49 | 32 | 13 | 36 |
| 30 | 35 | 50 | 15 |

C

| 6  | 59 | 42 | 23 |
|----|----|----|----|
| 57 | 24 | 5  | 44 |
| 12 | 37 | 56 | 25 |
| 55 | 10 | 27 | 38 |

D

Figure 2. Closed (or reentrant) magic tour of knight in 4 x 4 x 4 cube

| 19 | 46 | 63 | 2 |
|----|----|----|----|
| 48 | 17 | 4  | 61 |
| 29 | 52 | 33 | 16 |
| 34 | 15 | 30 | 51 |

A

| 58 | 7  | 22 | 43 |
|----|----|----|----|
| 5  | 60 | 41 | 24 |
| 40 | 9  | 28 | 53 |
| 27 | 54 | 39 | 10 |

B

| 47 | 18 | 3  | 62 |
|----|----|----|----|
| 20 | 45 | 64 | 1  |
| 49 | 32 | 13 | 36 |
| 14 | 35 | 50 | 31 |

C

| 6  | 59 | 42 | 23 |
|----|----|----|----|
| 57 | 8  | 21 | 44 |
| 12 | 37 | 56 | 25 |
| 55 | 26 | 11 | 38 |

D

Figure 3. Open magic tour of knight in 4 x 4 x 4 cube

## 3   Knight's tours in higher dimensions

This is a new topic in the age-old world of knight's tours. Manning [14] mentions that "The notion of geometries of n dimensions began to suggest itself to mathematicians about the middle of the 19$^{th}$ century. Cayley, Grassmann, Riemann, Clifford and some others introduced it into their mathematical investigations." Before we extend the knight's tour to four dimensions it is important for the reader to try to visualize a 4D hypercube. If we move a unit square a unit distance orthogonal to its plane in 3D space and join the corresponding corners, we get a cube. Analogously, we can imagine moving this unit cube a unit distance in an "orthogonal" direction in 4D space to produce the 4D equivalent of a 3D cube, which is known as a "hypercube". However we can only show this by means of a perspective drawing, as shown in Figure 4a.

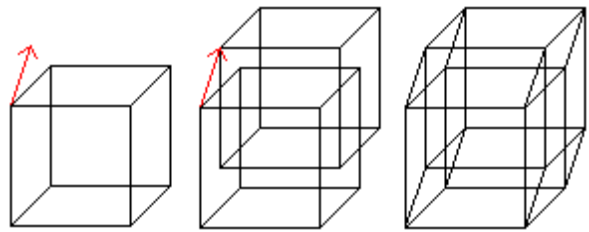

Figure 4a. Visualisation of a hypercube

This way, we can visualize a hypercube and be comfortable (and confident) with its 'look'. The hypercube has 16 corners (derived from 2 cubes), 32 edges (2 cubes and joining lines) and 24 square faces. Figure 4b shows all the possible knight's moves from a central cell in a $4 \times 4 \times 4 \times 4$ hypercube. Once the readers can visualize the jumps of the knight in hyperspace, they can count the possible number of knight moves from the nine distinct cell positions (lettered A to I) in the six planes (*xy*, *xz*, *xw*, *yz*, *yw*, *zw*) determined by pairs of the four coordinates *x*, *y*, *z*, *w*, as shown in Figure 5. On a two-dimensional board, knight has a choice of minimum 2 cells (when it is in the corner) and maximum up to 8 cells to which it can move. Now, readers can see that in 4D, knight has a choice of minimum 6x2 (= 12 cells) and maximum up to 6x8 (= 48 cells) to which it can move. In general, knight has a choice of minimum $n(n-1)$ cells and maximum up to $4n(n-1)$ cells to which it can move in a $n$-dimensional hypercube.

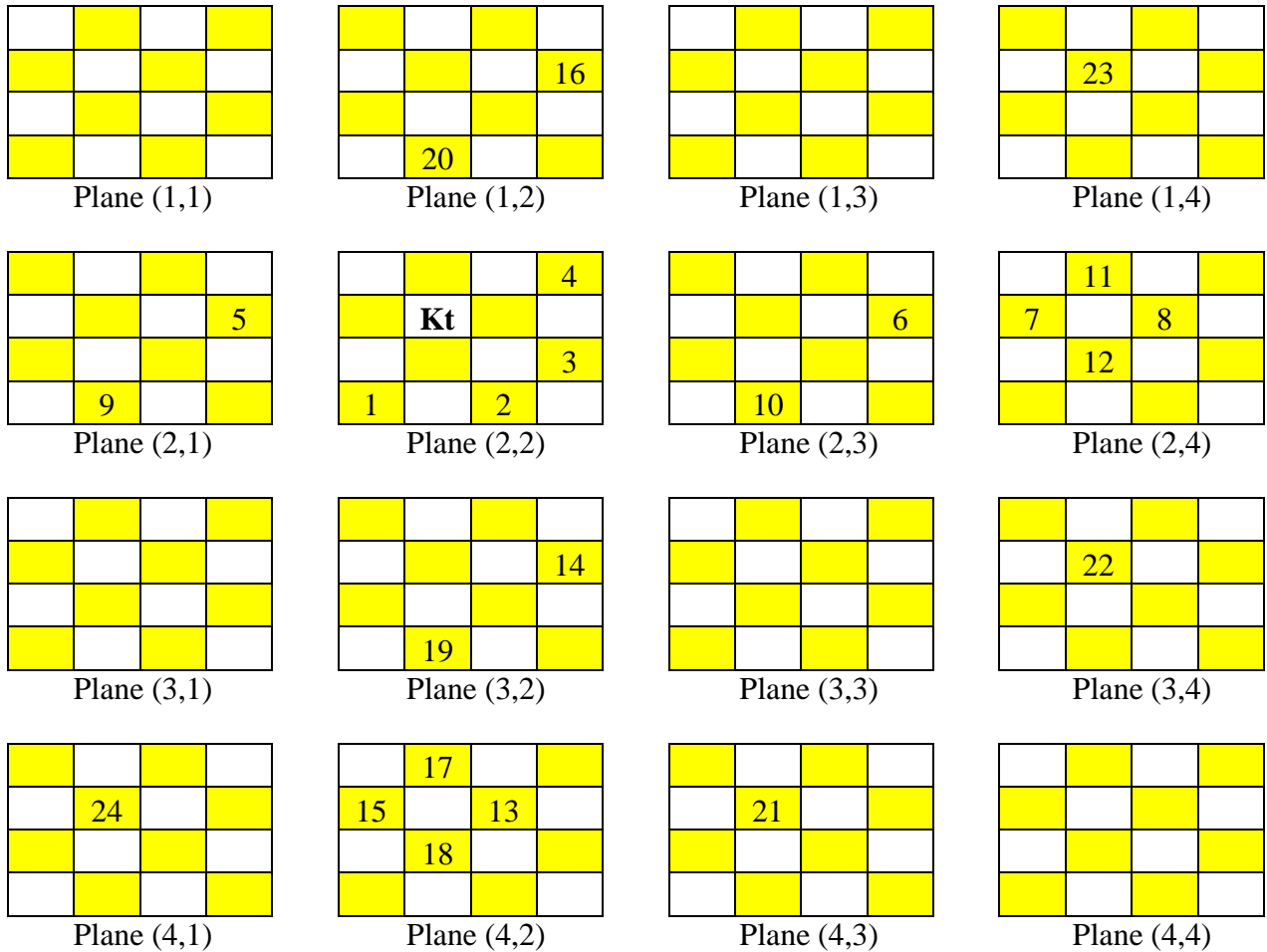

Figure 4b. Possible knight's move in 4 x 4 x 4 x 4 hypercube

After having a clear picture of the possible knight moves in hyperspace, we come to the tour of the knight in it. Readers can verify that $3 \times 4$ (= 12 cells) is the smallest rectangular lattice in which an open knight's tour is possible in 2-dimension. ($2 \times 3$ lattice is too small for the knight's tour and the knight can't enter or come out of the centre in a $3 \times 3$ lattice). Awani Kumar [15] has shown that $3 \times 4 \times 2$ (= 24 cells) is the smallest cuboid in which both closed and open tours are possible in three dimensions. Keeping these facts in mind, the author has observed that $3 \times 4 \times 2 \times 2$ (= 48 cells) is the smallest cuboid in which both closed and open tours are possible in four dimensions, as shown in Figure 6 and Figure 7 respectively. It can be further extended in 5-dimension cuboid of size $3 \times 4 \times 2 \times 2 \times 2$ (= 96 cells) as shown in Figure 8. In general, $3 \times 4 \times 2^{n-2}$ is the smallest cuboid in which knight's tour is possible in $n$-dimensions. Knight cannot move at all in a 2 x 2 x 2 x 2 cube and can neither move in nor come out from the central cell of a 3 x 3 x 3 x 3 cube. So, 4 x 4 x 4 x 4 is the smallest cube in which both closed and open tours are possible in 4-dimensions. In general, $4 \times 4 \times 4^{n-2}$ is the smallest cube in which knight's tour is possible in $n$-dimensions ($n \geq 3$).

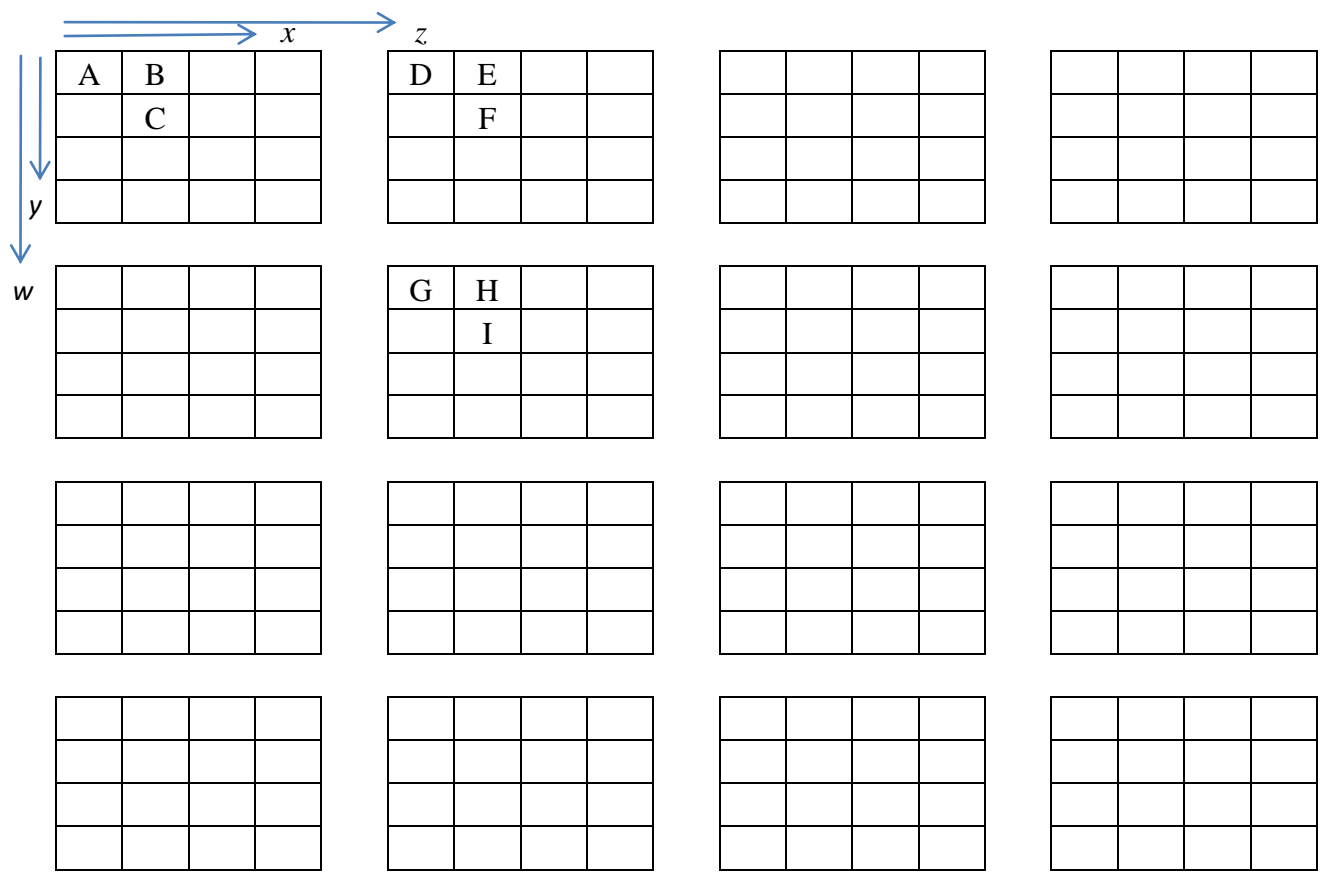

|       | A  | B  | C  | D  | E  | F  | G  | H  | I  |
|-------|----|----|----|----|----|----|----|----|----|
| $xy =$ | 2  | 3  | 4  | 2  | 3  | 4  | 2  | 3  | 4  |
| $xz =$ | 2  | 3  | 3  | 3  | 4  | 4  | 3  | 4  | 4  |
| $yz =$ | 2  | 2  | 3  | 3  | 3  | 4  | 3  | 3  | 4  |
| $xw =$ | 2  | 3  | 3  | 2  | 3  | 3  | 3  | 4  | 4  |
| $yw =$ | 2  | 2  | 3  | 2  | 2  | 3  | 3  | 3  | 4  |
| $zw =$ | 2  | 2  | 2  | 3  | 3  | 3  | 4  | 4  | 4  |
| Total | **12** | **15** | **18** | **15** | **18** | **21** | **18** | **21** | **24** |

Figure 5. Possible number of knight's moves from various cells in a 4 x 4 x 4 x 4 hypercube

| 1  | 22 | 9  | 18 |
|----|----|----|----|
| 8  | 17 | 2  | 23 |
| 3  | 24 | 7  | 16 |

| 10 | 13 | 6  | 21 |
|----|----|----|----|
| 5  | 20 | 11 | 14 |
| 12 | 15 | 4  | 19 |

| 46 | 31 | 42 | 27 |
|----|----|----|----|
| 41 | 26 | 47 | 32 |
| 48 | 33 | 40 | 25 |

| 43 | 28 | 39 | 36 |
|----|----|----|----|
| 38 | 35 | 44 | 29 |
| 45 | 30 | 37 | 34 |

Figure 6. Smallest knight's tour (closed) in a 3 x 4 x 2 x 2 hypercuboid

| 1  | 22 | 25 | 28 |
|----|----|----|----|
| 26 | 29 | 2  | 23 |
| 3  | 24 | 27 | 30 |

| 36 | 33 | 6  | 21 |
|----|----|----|----|
| 5  | 20 | 35 | 32 |
| 34 | 31 | 4  | 19 |

| 10 | 13 | 46 | 43 |
|----|----|----|----|
| 47 | 44 | 11 | 14 |
| 12 | 15 | 48 | 45 |

| 39 | 42 | 9  | 18 |
|----|----|----|----|
| 8  | 17 | 38 | 41 |
| 37 | 40 | 7  | 16 |

Figure 7. Smallest knight's tour (open) in a 3 x 4 x 2 x 2 hypercuboid

| 1  | 22 | 25 | 28 |
|----|----|----|----|
| 26 | 29 | 2  | 23 |
| 3  | 24 | 27 | 30 |

| 36 | 33 | 6  | 21 |
|----|----|----|----|
| 5  | 20 | 35 | 32 |
| 34 | 31 | 4  | 19 |

| 85 | 88 | 55 | 64 |
|----|----|----|----|
| 56 | 65 | 86 | 89 |
| 87 | 90 | 57 | 66 |

| 60 | 63 | 96 | 93 |
|----|----|----|----|
| 95 | 92 | 59 | 62 |
| 58 | 61 | 94 | 91 |

| 10 | 13 | 46 | 43 |
|----|----|----|----|
| 47 | 44 | 11 | 14 |
| 12 | 15 | 48 | 45 |

| 39 | 42 | 9  | 18 |
|----|----|----|----|
| 8  | 17 | 38 | 41 |
| 37 | 40 | 7  | 16 |

| 82 | 79 | 52 | 67 |
|----|----|----|----|
| 53 | 68 | 83 | 80 |
| 84 | 81 | 54 | 69 |

| 51 | 72 | 75 | 78 |
|----|----|----|----|
| 74 | 77 | 50 | 71 |
| 49 | 70 | 73 | 76 |

Figure 8. Smallest knight's tour (closed) in a five-dimensions hypercuboid

Schwenk [16], DeMaio and Mathew [17] have proved following theorems for closed knight's tour in two-dimensions and three-dimensions respectively.

**Theorem 1 (Schwenk)** *An $m \times n$ chessboard with $m \leq n$ has a closed knight's tour unless one or more of the following three conditions hold:*

(a) m and n are both odd;
(b) m ε {1, 2, 4};
(c) m = 3 and n ε {4, 6, 8}.

**Theorem 2 (DeMaio and Mathew)** *An $i \times j \times k$ chessboard for integers $i, j, k \geq 2$ has a closed knight's tour unless, without loss of generality, one or more of the following three conditions hold:*

(a) i, j and k are all odd;
(b) i = j = 2;
(c) i = 2 and j = k = 3.

Based on observations of knight's tour in higher dimensions, the author proposes to extend the result to four dimensional i x j x k x l cuboid for integers i, j, k, l ≥ 2 as follows:

**Conjecture 1 (Awani Kumar)** *An $i \times j \times k \times l$ chessboard for integers $i, j, k, l \geq 2$ has a closed knight's tour unless, one or more of the following three conditions hold:*

    (a) i, j, k and l are all odd;
    (b) i = j = k = 2;
    (c) i = j = 2 and k = l = 3.

The author proposes to further extend the result to five dimensional i x j x k x l x m cuboid for integers i, j, k, l, m ≥ 2 as follows:

**Conjecture 2 (Awani Kumar)** *An $i \times j \times k \times l \times m$ chessboard for integers $i, j, k, l, m \geq 2$ has a closed knight's tour unless, one or more of the following three conditions hold:*

    (a) i, j, k, l and m are all odd;
    (b) i = j = k = l = 2;
    (c) i = j = k = 2 and l = m = 3.

In general, the author proposes following results in an n-dimensional $a_1$ x $a_2$ x $a_3$ x … x $a_n$ cuboid for integers $a_1, a_2, a_3, \ldots, a_n \geq 2$ as follows:

**Conjecture 3 (Awani Kumar)** *An $a_1$ x $a_2$ x $a_3$ x … x $a_n$ chessboard for integers $a_1, a_2, a_3, \ldots, a_n \geq 2$ has a closed knight's tour unless, one or more of the following three conditions hold:*

    (a) $a_1, a_2, a_3, \ldots, a_{n-1}$ and $a_n$ are all odd;
    (b) $a_1 = a_2 = a_3 = \ldots = a_{n-1} = 2$;
    (c) $a_1 = a_2 = a_3 = \ldots = a_{n-2} = 2$ and $a_{n-1} = a_n = 3$.

## 4 Magic knight's tours in higher dimensions

There are trillions of knight's tours in the 4th order 4D hypercube and their number increases very rapidly with the order of the hypercube, therefore, they are not that challenging as such. However, magic knight's tours are a different story. It is more interesting (and challenging) to get magic knight's tours in 4 and higher dimensions. Before proceeding further, let us prove the following theorem.

**Theorem**: *There cannot be a magic knight's tour in a 4-dimensional hypercube of odd order.*
Proof: There are $m^2$ cells and $2m$ lines (rows and columns) in a square of order m. There are $m^3$ cells and $3m^2$ lines (rows, columns and pillars) in a cube of order m. Analogously, there are $m^4$ cells and $4m^3$ lines (rows, columns, pillars and posts) in a 4D hypercube of order m. Sum of the cells = $m^4(m^4 + 1)/2$. For the magic sum, the sum of each line = $m(m^4 + 1)/2$. Since m is odd, therefore, the magic sum will always be odd. We know that knight moves alternately from dark cell to light cell and vice-versa. So if the knight's tour starts from a dark cell then all the odd numbers will be in dark cells and all the even numbers will be in light cells. For odd m, if a line has (m − 1)/2 dark cells and (m + 1)/2 light cells then its neighboring line will have (m − 1)/2 light cells and (m + 1)/2 dark cells. Since these lines will have odd and even sum, so both can't sum to the magic constant. Proved.

Extending this logic, we see that there can't be a magic tour of knight in an *n*-dimensional hypercube of order *m* if *m* is odd. Similarly, there can't be a closed knight's tour in a 4-dimensional hypercube of order *m* if *m* is odd because the initial cell (1) and the last cell ($m^4$) will be of the same color. So they can't be connected by a knight's move. In general, there can't be a closed knight's tour in a *n*-dimensional hypercuboid having odd number of cells.

Can we have magic tour of knight in 4 × 4 × 4 × 4 hypercube? Yes!! As shown in Figure 9. But how can such tours be constructed? Well, the secret lies in the systematic movement of knight over various planes of the hypercube. Altogether, there are 256 cells. Put them into four groups: 1 to 64, 65 to 128, 129 to 192 and 193 to 256. Start from any cell and move the knight in such a way that each plane has equal numbers from each group in all the rows, columns, pillars and posts. With little patience and perseverance, magic tours can be constructed. The author got it by working on MS Excel. Here, sum of all the rows, columns, pillars and posts are 514. The author has enumerated over 200 such magic tours but in spite of intense effort, couldn't get the eight space diagonals magic and conjectures that such a tour doesn't exist. Readers having programming skills and access to powerful computers are requested to enumerate all such magic tours. Can it be extended in a 4D hypercube of order 8? Readers are encouraged to do so. What about magic tours in a hypercube of singly-even order (that is, 6, 10, 14 and so forth)? Well, perhaps, such magic tours don't exist. Jelliss [18] has proved that a magic knight's tour is impossible on a plane board with singly-even sides. But how close can we get to the magic tour? The best the author could get is 75% magic ratio on a 10 × 10 board, the highest on any singly-even board to date [19].

After constructing magic tour in a 4D hypercube, the author extended it into a 5-dimensional 4 × 4 × 4 × 4 × 4 hypercube as shown in Figure 10. Later, Nakamura [20] further extended it into 7-dimension hypercube but a magic tour in 6-dimensions is yet to be constructed.

## 5  Conclusion

Kaku [21] asserts that "There is a growing acknowledgment among physicists worldwide, including several Nobel laureates, that the universe may actually exist in higher-dimensional space." Pickover [22] declares "Various modern theories of "hyperspace" suggest that dimensions exist beyond the commonly accepted dimensions of space and time. The entire universe may actually exist in a higher-dimensional space. This idea is not science fiction …" Musser [23] muses that "As fantastic as extra dimensions of space sound, they might really exist… various mysteries of the world around us give the impression that the known universe is but the shadow of a higher-dimensional reality." Watkins [24] exclaims "… there is no reason to stop with chessboards in only three dimensions!" Accordingly, the author has shown possibility of knight's tour in higher dimensions and hopes that its study will help in unraveling secrets of hyperspace.

| 1 | 80 | 191 | 242 |
|---|---|---|---|
| 192 | 241 | 2 | 79 |
| 113 | 64 | 207 | 130 |
| 208 | 129 | 114 | 63 |

| 112 | 33 | 210 | 159 |
|---|---|---|---|
| 209 | 160 | 111 | 34 |
| 32 | 81 | 162 | 239 |
| 161 | 240 | 31 | 82 |

| 177 | 256 | 15 | 66 |
|---|---|---|---|
| 16 | 65 | 178 | 255 |
| 193 | 144 | 127 | 50 |
| 128 | 49 | 194 | 143 |

| 224 | 145 | 98 | 47 |
|---|---|---|---|
| 97 | 48 | 223 | 146 |
| 176 | 225 | 18 | 95 |
| 17 | 96 | 175 | 226 |

| 120 | 57 | 202 | 135 |
|---|---|---|---|
| 201 | 136 | 119 | 58 |
| 8 | 245 | 186 | 75 |
| 185 | 76 | 7 | 246 |

| 25 | 88 | 167 | 234 |
|---|---|---|---|
| 168 | 233 | 26 | 87 |
| 105 | 156 | 215 | 38 |
| 216 | 37 | 106 | 155 |

| 200 | 137 | 122 | 55 |
|---|---|---|---|
| 121 | 56 | 199 | 138 |
| 184 | 69 | 10 | 251 |
| 9 | 252 | 183 | 70 |

| 169 | 232 | 23 | 90 |
|---|---|---|---|
| 24 | 89 | 170 | 231 |
| 217 | 44 | 103 | 150 |
| 104 | 149 | 218 | 43 |

| 189 | 244 | 3 | 78 |
|---|---|---|---|
| 4 | 77 | 190 | 243 |
| 205 | 132 | 115 | 62 |
| 116 | 61 | 206 | 131 |

| 212 | 157 | 110 | 35 |
|---|---|---|---|
| 109 | 36 | 211 | 158 |
| 164 | 237 | 30 | 83 |
| 29 | 84 | 163 | 238 |

| 13 | 68 | 179 | 254 |
|---|---|---|---|
| 180 | 253 | 14 | 67 |
| 125 | 52 | 195 | 142 |
| 196 | 141 | 126 | 51 |

| 100 | 45 | 222 | 147 |
|---|---|---|---|
| 221 | 148 | 99 | 46 |
| 20 | 93 | 174 | 227 |
| 173 | 228 | 19 | 94 |

| 204 | 133 | 118 | 59 |
|---|---|---|---|
| 117 | 60 | 203 | 134 |
| 188 | 73 | 6 | 247 |
| 5 | 248 | 187 | 74 |

| 165 | 236 | 27 | 86 |
|---|---|---|---|
| 28 | 85 | 166 | 235 |
| 213 | 40 | 107 | 154 |
| 108 | 153 | 214 | 39 |

| 124 | 53 | 198 | 139 |
|---|---|---|---|
| 197 | 140 | 123 | 54 |
| 12 | 249 | 182 | 71 |
| 181 | 72 | 11 | 250 |

| 21 | 92 | 171 | 230 |
|---|---|---|---|
| 172 | 229 | 22 | 91 |
| 101 | 152 | 219 | 42 |
| 220 | 41 | 102 | 151 |

Figure 9. Magic tour of knight in a 4 x 4 x 4 x 4 hypercube

| 1 | 320 | 767 | 962 |
|---|---|---|---|
| 768 | 961 | 2 | 319 |
| 449 | 256 | 831 | 514 |
| 832 | 513 | 450 | 255 |

| 448 | 129 | 834 | 639 |
|---|---|---|---|
| 833 | 640 | 447 | 130 |
| 128 | 321 | 642 | 959 |
| 641 | 960 | 127 | 322 |

| 705 | 1024 | 63 | 258 |
|---|---|---|---|
| 64 | 257 | 706 | 1023 |
| 769 | 576 | 511 | 194 |
| 512 | 193 | 770 | 575 |

| 896 | 577 | 386 | 191 |
|---|---|---|---|
| 385 | 192 | 895 | 578 |
| 704 | 897 | 66 | 383 |
| 65 | 384 | 703 | 898 |

| 480 | 225 | 802 | 543 |
|---|---|---|---|
| 801 | 544 | 479 | 226 |
| 32 | 289 | 738 | 991 |
| 737 | 992 | 31 | 290 |

| 97 | 352 | 671 | 930 |
|---|---|---|---|
| 672 | 929 | 98 | 351 |
| 417 | 160 | 863 | 610 |
| 864 | 609 | 418 | 159 |

| 800 | 545 | 482 | 223 |
|---|---|---|---|
| 481 | 224 | 799 | 546 |
| 736 | 993 | 34 | 287 |
| 33 | 288 | 735 | 994 |

| 673 | 928 | 95 | 354 |
|---|---|---|---|
| 96 | 353 | 674 | 927 |
| 865 | 608 | 415 | 162 |
| 416 | 161 | 866 | 607 |

| 753 | 976 | 15 | 306 |
|---|---|---|---|
| 16 | 305 | 754 | 975 |
| 817 | 528 | 463 | 242 |
| 464 | 241 | 818 | 527 |

| 848 | 625 | 434 | 143 |
|---|---|---|---|
| 433 | 144 | 847 | 626 |
| 656 | 945 | 114 | 335 |
| 113 | 336 | 655 | 946 |

| 49 | 272 | 719 | 1010 |
|---|---|---|---|
| 720 | 1009 | 50 | 271 |
| 497 | 208 | 783 | 562 |
| 784 | 561 | 498 | 207 |

| 400 | 177 | 882 | 591 |
|---|---|---|---|
| 881 | 592 | 399 | 178 |
| 80 | 369 | 690 | 911 |
| 689 | 912 | 79 | 370 |

| 816 | 529 | 466 | 239 |
|---|---|---|---|
| 465 | 240 | 815 | 530 |
| 752 | 977 | 18 | 303 |
| 17 | 304 | 751 | 978 |

| 657 | 944 | 111 | 338 |
|---|---|---|---|
| 112 | 337 | 658 | 943 |
| 849 | 624 | 431 | 146 |
| 432 | 145 | 850 | 623 |

| 496 | 209 | 786 | 559 |
|---|---|---|---|
| 785 | 560 | 495 | 210 |
| 48 | 273 | 722 | 1007 |
| 721 | 1008 | 47 | 274 |

| 81 | 368 | 687 | 914 |
|---|---|---|---|
| 688 | 913 | 82 | 367 |
| 401 | 176 | 879 | 594 |
| 880 | 593 | 402 | 175 |

| 456 | 249 | 826 | 519 |
|---|---|---|---|
| 825 | 520 | 455 | 250 |
| 8 | 313 | 762 | 967 |
| 761 | 968 | 7 | 314 |

| 121 | 328 | 647 | 954 |
|---|---|---|---|
| 648 | 953 | 122 | 327 |
| 441 | 136 | 839 | 634 |
| 840 | 633 | 442 | 135 |

| 776 | 569 | 506 | 199 |
|---|---|---|---|
| 505 | 200 | 775 | 570 |
| 712 | 1017 | 58 | 263 |
| 57 | 264 | 711 | 1018 |

| 697 | 904 | 71 | 378 |
|---|---|---|---|
| 72 | 377 | 698 | 903 |
| 889 | 584 | 391 | 186 |
| 392 | 185 | 890 | 583 |

| 25 | 296 | 743 | 986 |
|---|---|---|---|
| 744 | 985 | 26 | 295 |
| 473 | 232 | 807 | 538 |
| 808 | 537 | 474 | 231 |

| 424 | 153 | 858 | 615 |
|---|---|---|---|
| 857 | 616 | 423 | 154 |
| 104 | 345 | 666 | 935 |
| 665 | 936 | 103 | 346 |

| 729 | 1000 | 39 | 282 |
|---|---|---|---|
| 40 | 281 | 730 | 999 |
| 793 | 552 | 487 | 218 |
| 488 | 217 | 794 | 551 |

| 872 | 601 | 410 | 167 |
|---|---|---|---|
| 409 | 168 | 871 | 602 |
| 680 | 921 | 90 | 359 |
| 89 | 360 | 679 | 922 |

| 824 | 521 | 458 | 247 |
|---|---|---|---|
| 457 | 248 | 823 | 522 |
| 760 | 969 | 10 | 311 |
| 9 | 312 | 759 | 970 |

| 649 | 952 | 119 | 330 |
|---|---|---|---|
| 120 | 329 | 650 | 951 |
| 841 | 632 | 439 | 138 |
| 440 | 137 | 842 | 631 |

| 504 | 201 | 778 | 567 |
|---|---|---|---|
| 777 | 568 | 503 | 202 |
| 56 | 265 | 714 | 1015 |
| 713 | 1016 | 55 | 266 |

| 73 | 376 | 695 | 906 |
|---|---|---|---|
| 696 | 905 | 74 | 375 |
| 393 | 184 | 887 | 586 |
| 888 | 585 | 394 | 183 |

| 745 | 984 | 23 | 298 |
|---|---|---|---|
| 24 | 297 | 746 | 983 |
| 809 | 536 | 471 | 234 |
| 472 | 233 | 810 | 535 |

| 856 | 617 | 426 | 151 |
|---|---|---|---|
| 425 | 152 | 855 | 618 |
| 664 | 937 | 106 | 343 |
| 105 | 344 | 663 | 938 |

| 41 | 280 | 727 | 1002 |
|---|---|---|---|
| 728 | 1001 | 42 | 279 |
| 489 | 216 | 791 | 554 |
| 792 | 553 | 490 | 215 |

| 408 | 169 | 874 | 599 |
|---|---|---|---|
| 873 | 600 | 407 | 170 |
| 88 | 361 | 682 | 919 |
| 681 | 920 | 87 | 362 |

| 765 | 964 | 3 | 318 |
|---|---|---|---|
| 4 | 317 | 766 | 963 |
| 829 | 516 | 451 | 254 |
| 452 | 253 | 830 | 515 |

| 836 | 637 | 446 | 131 |
|---|---|---|---|
| 445 | 132 | 835 | 638 |
| 644 | 957 | 126 | 323 |
| 125 | 324 | 643 | 958 |

| 61 | 260 | 707 | 1022 |
|---|---|---|---|
| 708 | 1021 | 62 | 259 |
| 509 | 196 | 771 | 574 |
| 772 | 573 | 510 | 195 |

| 388 | 189 | 894 | 579 |
|---|---|---|---|
| 893 | 580 | 387 | 190 |
| 68 | 381 | 702 | 899 |
| 701 | 900 | 67 | 382 |

| 804 | 541 | 478 | 227 |
|---|---|---|---|
| 477 | 228 | 803 | 542 |
| 740 | 989 | 30 | 291 |
| 29 | 292 | 739 | 990 |

| 669 | 932 | 99 | 350 |
|---|---|---|---|
| 100 | 349 | 670 | 931 |
| 861 | 612 | 419 | 158 |
| 420 | 157 | 862 | 611 |

| 484 | 221 | 798 | 547 |
|---|---|---|---|
| 797 | 548 | 483 | 222 |
| 36 | 285 | 734 | 995 |
| 733 | 996 | 35 | 286 |

| 93 | 356 | 675 | 926 |
|---|---|---|---|
| 676 | 925 | 94 | 355 |
| 413 | 164 | 867 | 606 |
| 868 | 605 | 414 | 163 |

| 13 | 308 | 755 | 974 |
|---|---|---|---|
| 756 | 973 | 14 | 307 |
| 461 | 244 | 819 | 526 |
| 820 | 525 | 462 | 243 |

| 436 | 141 | 846 | 627 |
|---|---|---|---|
| 845 | 628 | 435 | 142 |
| 116 | 333 | 654 | 947 |
| 653 | 948 | 115 | 334 |

| 717 | 1012 | 51 | 270 |
|---|---|---|---|
| 52 | 269 | 718 | 1011 |
| 781 | 564 | 499 | 206 |
| 500 | 205 | 782 | 563 |

| 884 | 589 | 398 | 179 |
|---|---|---|---|
| 397 | 180 | 883 | 590 |
| 692 | 909 | 78 | 371 |
| 77 | 372 | 691 | 910 |

| 468 | 237 | 814 | 531 |
|---|---|---|---|
| 813 | 532 | 467 | 238 |
| 20 | 301 | 750 | 979 |
| 749 | 980 | 19 | 302 |

| 109 | 340 | 659 | 942 |
|---|---|---|---|
| 660 | 941 | 110 | 339 |
| 429 | 148 | 851 | 622 |
| 852 | 621 | 430 | 147 |

| 788 | 557 | 494 | 211 |
|---|---|---|---|
| 493 | 212 | 787 | 558 |
| 724 | 1005 | 46 | 275 |
| 45 | 276 | 723 | 1006 |

| 685 | 916 | 83 | 366 |
|---|---|---|---|
| 84 | 365 | 686 | 915 |
| 877 | 596 | 403 | 174 |
| 404 | 173 | 878 | 595 |

| 828 | 517 | 454 | 251 |
|---|---|---|---|
| 453 | 252 | 827 | 518 |
| 764 | 965 | 6 | 315 |
| 5 | 316 | 763 | 966 |

| 645 | 956 | 123 | 326 |
|---|---|---|---|
| 124 | 325 | 646 | 955 |
| 837 | 636 | 443 | 134 |
| 444 | 133 | 838 | 635 |

| 508 | 197 | 774 | 571 |
|---|---|---|---|
| 773 | 572 | 507 | 198 |
| 60 | 261 | 710 | 1019 |
| 709 | 1020 | 59 | 262 |

| 69 | 380 | 699 | 902 |
|---|---|---|---|
| 700 | 901 | 70 | 379 |
| 389 | 188 | 891 | 582 |
| 892 | 581 | 390 | 187 |

| 741 | 988 | 27 | 294 |
|---|---|---|---|
| 28 | 293 | 742 | 987 |
| 805 | 540 | 475 | 230 |
| 476 | 229 | 806 | 539 |

| 860 | 613 | 422 | 155 |
|---|---|---|---|
| 421 | 156 | 859 | 614 |
| 668 | 933 | 102 | 347 |
| 101 | 348 | 667 | 934 |

| 37 | 284 | 731 | 998 |
|---|---|---|---|
| 732 | 997 | 38 | 283 |
| 485 | 220 | 795 | 550 |
| 796 | 549 | 486 | 219 |

| 412 | 165 | 870 | 603 |
|---|---|---|---|
| 869 | 604 | 411 | 166 |
| 92 | 357 | 678 | 923 |
| 677 | 924 | 91 | 358 |

| 460 | 245 | 822 | 523 |
|---|---|---|---|
| 821 | 524 | 459 | 246 |
| 12 | 309 | 758 | 971 |
| 757 | 972 | 11 | 310 |

| 117 | 332 | 651 | 950 |
|---|---|---|---|
| 652 | 949 | 118 | 331 |
| 437 | 140 | 843 | 630 |
| 844 | 629 | 438 | 139 |

| 780 | 565 | 502 | 203 |
|---|---|---|---|
| 501 | 204 | 779 | 566 |
| 716 | 1013 | 54 | 267 |
| 53 | 268 | 715 | 1014 |

| 693 | 908 | 75 | 374 |
|---|---|---|---|
| 76 | 373 | 694 | 907 |
| 885 | 588 | 395 | 182 |
| 396 | 181 | 886 | 587 |

| 21 | 300 | 747 | 982 |
|---|---|---|---|
| 748 | 981 | 22 | 299 |
| 469 | 236 | 811 | 534 |
| 812 | 533 | 470 | 235 |

| 428 | 149 | 854 | 619 |
|---|---|---|---|
| 853 | 620 | 427 | 150 |
| 108 | 341 | 662 | 939 |
| 661 | 940 | 107 | 342 |

| 725 | 1004 | 43 | 278 |
|---|---|---|---|
| 44 | 277 | 726 | 1003 |
| 789 | 556 | 491 | 214 |
| 492 | 213 | 790 | 555 |

| 876 | 597 | 406 | 171 |
|---|---|---|---|
| 405 | 172 | 875 | 598 |
| 684 | 917 | 86 | 363 |
| 85 | 364 | 683 | 918 |

Figure 10. Magic tour of knight in a 4 x 4 x 4 x 4 x 4 hypercube

## Acknowledgement

The author is grateful to Manoj Kumar Agarwal, Atulit Kumar and George Jelliss for their help in preparation of this article.

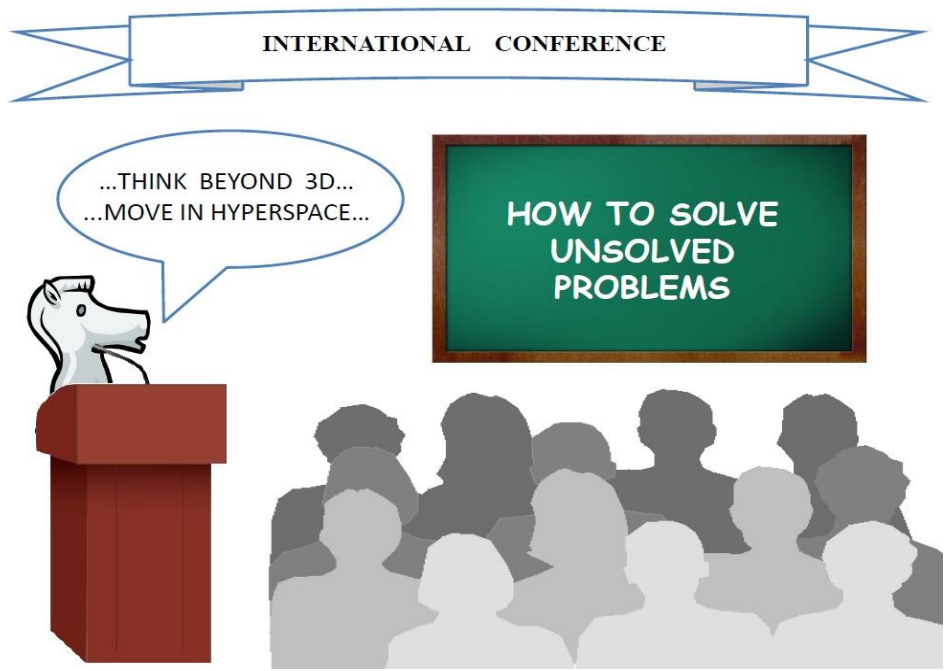